\documentclass{amsart}
\usepackage{amsmath}
\usepackage{bbm}
\usepackage{mathrsfs}
\usepackage{amssymb,amsmath,amsfonts,amsthm,graphicx,enumerate,amscd,latexsym,curves,multibox}
\usepackage[mathscr]{eucal}
\usepackage{epsfig,epsf,xypic,epic}

\usepackage[T1]{fontenc}
\usepackage[sc]{mathpazo}

\newtheorem{thm}{Theorem}[section]
\newtheorem{lem}{Lemma}[section]
\newtheorem{cor}{Corollary}[section]
\newtheorem{prop}{Proposition}[section]
\newtheorem{defn}{Definition}[section]
\newtheorem{nb}{Remark}[section]

\numberwithin{equation}{section}

\begin{document}

\title[Holomorphic line bundles from
Lagrangian sections by SYZ transformations]{Holomorphic line bundles
on projective toric manifolds from Lagrangian sections of their
mirrors by SYZ transformations}
\author[K.-W. Chan]{Kwokwai Chan}
\address{Department of Mathematics, Harvard University,
1 Oxford Street, Cambridge, MA 02138, USA}
\email{kwchan@math.harvard.edu}

\begin{abstract}
The mirror of a projective toric manifold $X_\Sigma$ is given by a
Landau-Ginzburg model $(Y,W)$. We introduce a class of Lagrangian
submanifolds in $(Y,W)$ and show that, under the SYZ mirror
transformation, they can be transformed to torus-invariant hermitian
metrics on holomorphic line bundles over $X_\Sigma$. Through this
geometric correspondence, we also identify the mirrors of
Hermitian-Einstein metrics, which are given by distinguished
Lagrangian sections whose potentials satisfy certain Laplace-type
equations.
\end{abstract}

\maketitle

\tableofcontents

\section{Introduction}

Let $X_\Sigma$ be a projective toric manifold defined by a fan
$\Sigma$. The mirror of $X_\Sigma$ is given by a Landau-Ginzburg
model $(Y,W)$, which consists of a noncompact K\"{a}hler manifold
$Y$ and a holomorphic function $W:Y\rightarrow\mathbb{C}$ (the
superpotential). Mirror symmetry relates the complex geometry of
$X_\Sigma$ to the symplectic geometry of $(Y,W)$. In particular,
holomorphic vector bundles (or more generally, coherent sheaves)
over $X_\Sigma$ should correspond to Lagrangian cycles in $(Y,W)$.
This is succinctly expressed by Kontsevich's Homological Mirror
Symmetry Conjecture for toric manifolds \cite{Kontsevich98}, which
states that the derived category of coherent sheaves
$D^bCoh(X_\Sigma)$ is equivalent to the Fukaya-Kontsevich-Seidel
category of $(Y,W)$. Since then, much work has been done
\cite{HIV00}, \cite{Seidel00}, \cite{Ueda04}, \cite{AKO04},
\cite{AKO05}, \cite{Abouzaid05}, \cite{Fang08}, culminating in
proofs of the conjecture for all projective toric manifolds in
Abouzaid \cite{Abouzaid06} and, more recently, in
Fang-Liu-Treumann-Zaslow
\cite{FLTZ08}.\footnote{Fang-Liu-Treumann-Zaslow \cite{FLTZ08} also
proved an equivariant version of the conjecture.}

In this paper, we will examine the correspondence between
holomorphic line bundles on $X_\Sigma$ and Lagrangian cycles on
$(Y,W)$ from a different angle, namely, by applying \textit{SYZ
mirror transformations} \cite{Chan-Leung08a}, \cite{Chan-Leung08b}.
Our goal is to put the correspondence in the toric case in the same
footing as the semi-flat Calabi-Yau case as done in Leung-Yau-Zaslow
\cite{LYZ00}. This approach is also closely related to the works
\cite{Abouzaid05}, \cite{Abouzaid06}, \cite{Fang08}, \cite{FLTZ08},
where T-duality was used implicitly or explicitly.

Let $N\cong\mathbb{Z}^n$ be a rank $n$ lattice,
$M=\mbox{Hom}(N,\mathbb{Z})$ the dual lattice and
$\langle\cdot,\cdot\rangle:M\times N\rightarrow\mathbb{Z}$ the dual
pairing, and let $N_\mathbb{R}=N\otimes_\mathbb{Z}\mathbb{R}$,
$M_\mathbb{R}=M\otimes_\mathbb{Z}\mathbb{R}$. Denote by $T_N$ and
$T_M$ the real tori $N_\mathbb{R}/N$ and $M_\mathbb{R}/M$
respectively. A projective toric $n$-fold $X_\Sigma$ contains an
open dense torus orbit
$U=N\otimes_\mathbb{Z}\mathbb{C}^*\cong(\mathbb{C}^*)^n$, which can
also be written as
$$U=N_\mathbb{R}\times\sqrt{-1}T_N=TN_\mathbb{R}/N,$$
where we have, by abuse of notations, also used $N$ to denote the
family of lattices $N_\mathbb{R}\times\sqrt{-1}N\subset
TN_\mathbb{R}$. The projection map $U\rightarrow N_\mathbb{R}$ is a
(trivial) torus bundle. According to the philosophy of the
\textit{Strominger-Yau-Zaslow Conjecture} \cite{SYZ96}, the mirror
manifold $Y$ is given by the dual torus bundle (see
\cite{Chan-Leung08a}, \cite{Chan-Leung08b})
$$Y=N_\mathbb{R}\times\sqrt{-1}T_M=T^*N_\mathbb{R}/M,$$
with $M$ denoting the family of lattices
$N_\mathbb{R}\times\sqrt{-1}M\subset T^*N_\mathbb{R}$. Using the
semi-flat SYZ mirror transformation (or T-duality), $T_N$-invariant
hermitian metrics on holomorphic line bundles over $X_\Sigma$ (when
restricted to $U$) can be transformed to give Lagrangian sections of
$Y\rightarrow N_\mathbb{R}$ as in \cite{LYZ00}.\footnote{More
precisely, one should get Lagrangian sections equipped with flat
$U(1)$-connections. But our Lagrangian sections are simply
connected, so all flat $U(1)$-connections are gauge equivalent to
the trivial one and we will ignore this data.} Naturally, one would
ask the following\\

\noindent\textbf{Question: }\textit{Which Lagrangian sections of
$Y\rightarrow N_\mathbb{R}$ can be transformed back, by the inverse
SYZ mirror transformation, to $T_N$-invariant hermitian metrics on
holomorphic line bundles over $X_\Sigma$?}\\

\noindent Put it in another way, the problem is to characterize the
set of Lagrangian sections of $Y\rightarrow N_\mathbb{R}$ we get by
transforming $T_N$-invariant hermitian metrics on holomorphic line
bundles over $X_\Sigma$. One of our aims in this paper is to answer
this question.

Recall that the superpotential $W$ is a Laurent polynomial (see, for
example, \cite{Chan-Leung08a}, \cite{Chan-Leung08b}). Write $W$ as a
sum of monomials: $W=\sum_{i=1}^d W_i$. In a sense, the monomial
$W_i$ (for $i=1,\ldots,d$) is mirror to the toric prime divisor
$D_i\subset\bar{X}$ associated to the primitive generator $v_i\in N$
of a 1-dimensional cone in $\Sigma$. Consider the embedding
$\iota:M\hookrightarrow\mathbb{Z}^d$ defined by $\iota(u)=(\langle
u,v_1\rangle,\ldots,\langle u,v_d\rangle)$. By the theory of toric
varieties, the quotient $\mathbb{Z}^d/\iota(M)$ is canonically
identified with $H^2(X_\Sigma,\mathbb{Z})$. In Section~\ref{sec3},
we will define, for each $[a]\in H^2(X_\Sigma,\mathbb{Z})$, a growth
condition $(\ast_{[a]})$ for Lagrangian sections of $Y\rightarrow
N_\mathbb{R}$. We can now state our main result as follows, which
will be proved in Section~\ref{sec4}.
\begin{thm}
Let $\mathcal{L}_{[a]}$ be the holomorphic line bundle over
$X_\Sigma$ corresponding to $[a]\in H^2(X_\Sigma,\mathbb{Z})$. Then
the SYZ mirror transformation gives a bijective correspondence
between $T_N$-invariant hermitian metrics on $\mathcal{L}_{[a]}$ and
Lagrangian sections of $Y\rightarrow N_\mathbb{R}$ satisfying the
growth condition $(\ast_{[a]})$.
\end{thm}
Notice that all Lagrangian sections of $Y\rightarrow N_\mathbb{R}$
are Hamiltonian isotopic to the zero section, i.e. they represent
the same Hamiltonian class. To get a correspondence with the class
of holomorphic line bundles on $X_\Sigma$, it is therefore necessary
to find a finer equivalence relation. For this purpose, we define
two Lagrangian sections of $(Y,W)$ to be equivalent if they can be
deformed to each other through Hamiltonian isotopies which
\textit{preserve a growth condition $(\ast_{[a]})$}. It is easy to
see that each equivalence class then consists of exactly those
Lagrangian sections which satisfy the same growth condition
$(\ast_{[a]})$.

Furthermore, by our main result, we can easily identify the
Lagrangian sections which are \textit{mirror to Hermitian-Einstein
metrics} on holomorphic line bundles. These turn out to be
Lagrangian sections whose potentials satisfy certain Laplace-type
equations. We call these Lagrangian sections \textit{harmonic}.
Hence, as an immediate consequence of our main result, we have the
following
\begin{cor}\mbox{}
\begin{enumerate}
\item[1.] The SYZ mirror transformation provides a bijective
correspondence between isomorphism classes of holomorphic line
bundles over $X_\Sigma$ and equivalence classes of Lagrangian
sections of $(Y,W)$.
\item[2.] Each equivalence class of Lagrangian sections of $(Y,W)$ is represented by
a unique harmonic Lagrangian section.
\end{enumerate}
\end{cor}
All of these will be discussed with more details in
Section~\ref{sec4}. The next section (Section~\ref{sec2}) is a brief
review of mirror symmetry for toric manifolds. Some further remarks
and discussions are contained in the final section
(Section~\ref{sec5}).\\

\noindent\textbf{Acknowledgments.} I am grateful to Siu-Cheong Lau
for numerous useful discussions. Comments from an anonymous referee
were very helpful and led to a significant improvement in the
exposition. I would also like to thank Professor Shing-Tung Yau and
Profesor Naichung Conan Leung for their continuous encouragement and
support. This work was supported by Harvard University and the
Croucher Foundation Fellowship.

\section{Projective toric manifolds and their mirrors}\label{sec2}

In this section, we briefly review the geometric aspects of the
mirror symmetry for projective toric manifolds and fix our
notations.

A projective toric manifold by $X_\Sigma$ is defined by a smooth,
complete fan $\Sigma$ in $N_\mathbb{R}$. By the general theory of
toric varieties \cite{Fulton93}, \cite{Guillemin94a}, any ample line
bundle $\mathcal{L}$ on $X_\Sigma$ is determined by a lattice
polytope $\bar{P}\subset M_\mathbb{R}$ dual to $\Sigma$. If
$v_1,\ldots,v_d\in N$ are the primitive generators of the
1-dimensional cones of $\Sigma$, then there is a $d$-tuple of
integers $\lambda=(\lambda_1,\ldots,\lambda_d)\in\mathbb{Z}^d$ such
that
$$\bar{P}=\{x=(x_1,\ldots,x_n)\in M_\mathbb{R}:
\langle x,v_i\rangle+\lambda_i\geq0\textrm{ for }i=1,\ldots,d\},$$
and  $\mathcal{L}$ is then canonically identified with the divisor
line bundle $\mathcal{O}(D_\lambda)$, where
$D_\lambda=\sum_{i=1}^d\lambda_iD_i$ is an ample toric divisor. We
fix such an ample line bundle $\mathcal{L}$ and equip $X_\Sigma$
with the K\"{a}hler structure
$\omega_{X_\Sigma}=\iota^*\omega_{FS}$, where
$\iota:X_\Sigma\hookrightarrow\mathbb{C}P^N$ is an embedding induced
by $\mathcal{L}$ (note that since $X_\Sigma$ is smooth and
projective, every ample line bundle $\mathcal{L}$ is in fact very
ample; see Fulton \cite{Fulton93}), and $\omega_{FS}$ is the
Fubini-Study K\"{a}hler structure on $\mathbb{C}P^N$.

Recall that $X_\Sigma$ contains an open dense orbit
$U=X_\Sigma\setminus\bigcup_{i=1}^d
D_i=N\otimes_\mathbb{Z}\mathbb{C}^*=N_\mathbb{R}\times\sqrt{-1}T_N=TN_\mathbb{R}/N$,
and we have a natural torus fibration
$\nu_U:U=TN_\mathbb{R}/N\rightarrow N_\mathbb{R}$ given by
projection to the first factor. If $\xi_1,\ldots,\xi_n\in\mathbb{R}$
and $u_1,\ldots,u_n\in\mathbb{R}/2\pi\mathbb{Z}$ are the base
coordinates on $N_\mathbb{R}$ and fiber coordinates on $T_N$
respectively, then the complex coordinates on $U=(\mathbb{C}^*)^n$
are given by $w_j=e^{\xi_j+\sqrt{-1}u_j}$, $j=1,\ldots,n$, and the
restriction of $\omega_{X_\Sigma}$ to $U$ can be explicitly written
as
$$\omega_U=\omega_{X_\Sigma}|_U=2\sqrt{-1}\partial\bar{\partial}\phi=
\sum_{j,k=1}^n\frac{\partial^2\phi}{\partial\xi_j\partial\xi_k}d\xi_j\wedge
du_k,$$ where $\phi:N_\mathbb{R}\rightarrow\mathbb{R}$ is the
function given by
$$\phi(\xi)=\frac{1}{2}\log\Bigg(\sum_{u\in\bar{P}\cap M}c_ue^{2\langle u,\xi\rangle}\Bigg),$$
for some nonnegative constants $c_u$, $u\in\bar{P}\cap M$, which
depend on the embedding $\iota$. We use $\phi_j$ and $\phi_{jk}$ to
denote the partial derivatives $\frac{\partial\phi}{\partial\xi_j}$
and $\frac{\partial^2\phi}{\partial\xi_j\partial\xi_k}$
respectively, and let $(\phi^{jk})_{j,k=1}^n$ be the inverse matrix
of $(\phi_{jk})_{j,k=1}^n$.

If $\mu:X_\Sigma\rightarrow\bar{P}$ is the moment map of the
Hamiltonian $T_N$-action on $(X_\Sigma,\omega_{X_\Sigma})$, then the
restriction of $\mu$ to $U\subset X_\Sigma$ is the map
$\mu_U:U\rightarrow M_\mathbb{R}$ given by
$$\mu_U(w)=d\phi(\log|w_1|,\ldots,\log|w_n|)=
\frac{\sum_{u\in\bar{P}\cap M}c_u|w^u|^2\cdot
u}{\sum_{u\in\bar{P}\cap M}c_u|w^u|^2},$$ for $w=(w_1,\ldots,w_n)\in
U=(\mathbb{C}^*)^n$. The image of $\mu_U$ is the interior $P$ of the
polytope $\bar{P}$. In fact, the Legendre transform of the function
$\phi$ gives a diffeomorphism $\Phi=d\phi:N_\mathbb{R}\rightarrow P$
and $\mu_U=\Phi\circ\nu_U$. We also have a nowhere vanishing
holomorphic $n$-form on $U$ given by
$$\Omega_U=\frac{dw_1}{w_1}\wedge\ldots\wedge\frac{dw_n}{w_n}.$$
With respect to $\omega_U$ and $\Omega_U$, $\nu_U:U\rightarrow
N_\mathbb{R}$ and $\mu_U:U\rightarrow P$ are special Lagrangian
torus fibrations, in the sense of Auroux \cite{Auroux07} (and $U$ is
an almost Calabi-Yau manifold).\\

The mirror of $X_\Sigma$ is the Landau-Ginzburg model $(Y,W)$
described as follows. The mirror manifold $Y$ is the dual torus
fibration $Y=N_\mathbb{R}\times\sqrt{-1}T_M=T^*N_\mathbb{R}/M$.
Written in this way, $Y$ is naturally a symplectic manifold,
equipped with the standard symplectic structure
$\omega_Y=\sum_{j=1}^nd\xi_j\wedge dy_j$, where
$y_1,\ldots,y_n\in\mathbb{R}/2\pi\mathbb{Z}$ are the dual
coordinates on the fiber $T_M$. The projection map
$\mu_Y:Y=T^*N_\mathbb{R}/M\rightarrow N_\mathbb{R}$ is the moment
map for the Hamiltonian $T_M$-action on $Y$. To describe the complex
structure on $Y$ and write down the superpotential $W$, it is more
convenient to change the coordinates on the base by the
diffeomorphism $\Phi:N_\mathbb{R}\rightarrow P$ and rewrite $Y$ as
$Y=P\times\sqrt{-1}T_M=TP/M$, where $M$ here denotes the (trivial)
family of lattices $P\times\sqrt{-1}M$. Then $Y$ is naturally a
complex manifold with complex coordinates given by
$z_j=e^{-x_j+\sqrt{-1}y_j}$, where $x_1,\ldots,x_n$ are the
coordinates on $P$. There is a nowhere vanishing holomorphic
$n$-form on $Y$ given by
$$\Omega_Y=\frac{dz_1}{z_1}\wedge\ldots\wedge\frac{dz_n}{z_n}.$$
The superpotential $W:Y\rightarrow\mathbb{C}$ is the Laurent
polynomial
$$W(z)=e^{-\lambda_1}z^{v_1}+\ldots+e^{-\lambda_d}z^{v_d}$$
for $z=(z_1,\ldots,z_n)\in(\mathbb{C}^*)^n$, where $z^{v_i}$ denotes
the monomial $z_1^{v_i^1}\ldots z_n^{v_i^n}$. $W$ can be obtained as
the SYZ mirror transformation of a certain function on the geodesic
loop space $L_U$ of $U\subset X_\Sigma$ (see Chan-Leung
\cite{Chan-Leung08a}, \cite{Chan-Leung08b} for details).

Notice that, as a complex manifold, $Y$ is biholomorphic to the
bounded domain
$\{z\in(\mathbb{C}^*)^n:|e^{-\lambda_i}z^{v_i}|<1,\textrm{ for
}i=1,\ldots,d\}$ in $(\mathbb{C}^*)^n$. On the other hand, since
$\Phi$ is a Legendre transform, there exists a function
$\psi:P\rightarrow\mathbb{R}$ such that
$\phi^{jk}=\psi_{jk}=\frac{\partial^2\psi}{\partial x_j\partial
x_k}$ and $(\psi^{jk}):=(\psi_{jk})^{-1}=(\phi_{jk})$. The Legendre
transform $\Psi:P\rightarrow N_\mathbb{R}$ of $\psi$ is then the
inverse of $\Phi:N_\mathbb{R}\rightarrow P$, i.e. $\Psi=\Phi^{-1}$.
Now, the symplectic structure $\omega_Y$ is given in the $x_j,y_j$
coordinates by
$$\omega_Y=\sum_{j,k=1}^n\frac{\partial^2\psi}{\partial x_j\partial
x_k}dx_j\wedge dy_k.$$ If we denote by $\nu_Y:Y=TP/M\rightarrow P$
the projection map to the base $P$, then we have
$\mu_Y=\Psi\circ\nu_Y$. With respect to $\omega_Y$ and $\Omega_Y$,
$\nu_Y:Y\rightarrow N_\mathbb{R}$ and $\mu_Y:Y\rightarrow P$ are
special Lagrangian torus fibrations, which are dual to
$\nu_U:U\rightarrow N_\mathbb{R}$ and $\mu_U:U\rightarrow P$
respectively.\\

Physical arguments predict that the complex (respectively,
symplectic) geometry of $X_\Sigma$ is interchanged with the
symplectic (respectively, complex) geometry of $(Y,W)$ under mirror
symmetry. For precise mathematical statements and how SYZ mirror
transformations are applied to explain the geometry underlying this
mirror symmetry, we refer the reader to \cite{Chan-Leung08a},
\cite{Chan-Leung08b}.

\section{A class of Lagrangian submanifolds in Landau-Ginzburg
models}\label{sec3}

In this section, we introduce a class of Lagrangian submanifolds in
$(Y,W)$, which are sections of the torus fibration
$\mu_Y:Y\rightarrow N_\mathbb{R}$ (or $\nu_Y:Y\rightarrow P$),
satisfying certain growth conditions at infinity.

Let $(Y,W)$ be a Landau-Ginzburg model mirror to a projective toric
manifold $X_\Sigma$. Recall that the superpotential
$W\in\mathcal{O}(Y)$ is a Laurent polynomial of the form
$\sum_{i=1}^d b_iz^{v_i}$, for some $v_1,\ldots,v_d\in N$. Define
$A(W)$ to be the quotient group $\mathbb{Z}^d/\iota(M)$, where
$\iota:M\hookrightarrow\mathbb{Z}^d$, $u\mapsto(\langle
u,v_1\rangle,\ldots,\langle u,v_d\rangle)$ is the homomorphism
defined in the introduction. As we have mentioned before, $A(W)$ is
canonically identified with the second cohomology group
$H^2(X_\Sigma,\mathbb{Z})$ of $X_\Sigma$. Moreover, if we let
$\textrm{Log}:TM_\mathbb{R}/M=(\mathbb{C}^*)^n\rightarrow
M_\mathbb{R}=\mathbb{R}^n$ be the map defined by
$$\textrm{Log}(z_1,\ldots,z_n)=(\log|z_1|,\ldots,\log|z_n|),$$
then (the closure of) the image of $Y$ under $\textrm{Log}$, i.e.
$P:=\textrm{Log}(Y)=\textrm{Log}(\{z\in(\mathbb{C}^*)^n:|b_vz^v|<1\textrm{
for all }v\in A\})$, is a polytope in $M_\mathbb{R}$, and this
determines a fan $\Sigma$ in $N_\mathbb{R}$. These are exactly the
polytope and fan defining the projective toric manifold $X_\Sigma$.

Now, we write $Y=N_\mathbb{R}\times\sqrt{-1}T_M=T^*N_\mathbb{R}/M$
and equip $Y$ with the standard symplectic form
$\omega_Y=\sum_{j=1}^nd\xi_j\wedge dy_j$. Since $N_\mathbb{R}$ is
simply connected, any section $L$ of $\mu_Y:Y\rightarrow
N_\mathbb{R}$ can be lifted to a section
$\tilde{L}=\{(\xi,y(\xi)):\xi=(\xi_1,\ldots,\xi_n)\in
N_\mathbb{R}\}$ of $T^*N_\mathbb{R}$, where
$y:N_\mathbb{R}\rightarrow M_\mathbb{R}$ should be regarded as a
1-form on $N_\mathbb{R}$; moreover, if $\{(\xi,y_1(\xi)):\xi\in
N_\mathbb{R}\}, \{(\xi,y_2(\xi)):\xi\in N_\mathbb{R}\}\subset
T^*N_\mathbb{R}$ are two lifts of $L\subset Y$, then $y_1-y_2\equiv
u$ for some constant $u\in M$. By the standard argument as shown in
\cite{LYZ00}, a section $L$ of $\mu_Y:Y\rightarrow N_\mathbb{R}$ is
Lagrangian if and only if some lift
$\tilde{L}=\{(\xi,y(\xi)):\xi=(\xi_1,\ldots,\xi_n)\in
N_\mathbb{R}\}$ of $L$ to $T^*N_\mathbb{R}$ is the graph of an exact
1-form, i.e. if and only if
$$y(\xi)=dg(\xi)
=\Bigg(\frac{\partial g}{\partial\xi_1},\ldots,\frac{\partial
g}{\partial\xi_n}\Bigg),$$ for some function $g$ on $N_\mathbb{R}$,
which is unique up to adding a constant. $g$ is called a
\textit{potential} of the lift $\tilde{L}$ of the Lagrangian section
$L$. For our purpose, we need $g$ to be of class $C^2$.
\begin{defn}\label{def3.1}
Let $a=(a_1,\ldots,a_d)\in\mathbb{Z}^d$ be a $d$-tuple of integers.
A Lagrangian section $\tilde{L}=\{(\xi,y(\xi)):\xi\in
N_\mathbb{R}\}$ of $T^*N_\mathbb{R}\rightarrow N_\mathbb{R}$ is said
to satisfy the growth condition $(\ast_a)$ if a potential $g\in
C^2(N_\mathbb{R})$ of $\tilde{L}$ satisfies the following
conditions: Given any $n$-dimensional cone $\sigma\in\Sigma$.
Suppose that, without loss of generality, $\sigma$ is generated by
$v_1,\ldots,v_n$; and let
$\xi(t)=\xi(t_1,\ldots,t_n)=t_1v_1+\ldots+t_nv_n$, for
$t=(t_1,\ldots,t_n)\in\mathbb{R}^n$. Then, we have,
\begin{enumerate}
\item[1.] the functions
$2e^{-2t_j}(\langle dg(\xi(t)),v_j\rangle+a_j)$ and
$e^{-2t_j}(v_j^T\textrm{Hess}(g)v_j)(\xi(t))$ have the same limit as
$t_j\rightarrow-\infty$, for $j=1,\ldots,n$;
\item[2.] for any $j,k,l\in\{1,\ldots,n\}$, the function
$(v_j^T\textrm{Hess}(g)v_k)(\xi(t))$ has a limit as
$t_l\rightarrow-\infty$; and,
\item[3.] for any distinct $j,k\in\{1,\ldots,n\}$, the function
$e^{-t_j-t_k}(v_j^T\textrm{Hess}(g)v_k)(\xi(t))$ goes to zero when
$t_j\rightarrow-\infty$ or $t_k\rightarrow-\infty$.
\end{enumerate}
Let $[a]\in A(W)$. A Lagrangian section $L$ of $\mu_Y:Y\rightarrow
N_\mathbb{R}$ is said to satisfy the growth condition $(\ast_{[a]})$
if some lift $\tilde{L}$ of $L$ from $Y$ to $T^*N_\mathbb{R}$
satisfies $(\ast_a)$ for some representative
$a=(a_1,\ldots,a_d)\in\mathbb{Z}^d$ of $[a]$.
\end{defn}
We denote the set of Lagrangian sections of $\mu_Y:Y\rightarrow
N_\mathbb{R}$ satisfying $(\ast_{[a]})$ for some $[a]\in A(W)$ by
$\mathbbm{L}(Y,W)$.
\begin{nb}
The condition that a Lagrangian section $L$ of $\mu_Y:Y\rightarrow
N_\mathbb{R}$ satisfies $(\ast_{[a]})$ is well-defined because if
$\{(\xi,y_1(\xi)):\xi\in N_\mathbb{R}\}, \{(\xi,y_2(\xi)):\xi\in
N_\mathbb{R}\}\subset T^*N_\mathbb{R}$ are two lifts of $L\subset
Y$, then their potentials $g_1,g_2$ will differ by a linear function
of the form $\langle u,\xi\rangle+\alpha$, for some $u\in M$ and
$\alpha\in\mathbb{R}$. Thus, when one of the lifts satisfies
$(\ast_a)$, the other will satisfy $(\ast_{a'})$, where
$a'=a+(\langle u,v_1\rangle,\ldots,\langle u,v_d\rangle)$, and note
that we have $[a]=[a']$.
\end{nb}
We give a couple of examples to illustrate our definitions.\\

\noindent\textbf{Example 1.} The simplest example is given by
$X_\Sigma=\mathbb{C}P^1$. The fan $\Sigma$ in
$N_\mathbb{R}=\mathbb{R}$ is generated by two primitive vectors
$v_1=1,v_2=-1$ (see Figure 1 below).
\begin{figure}[ht]
\setlength{\unitlength}{1mm}
\begin{picture}(100,8)
\put(50,4){\vector(1,0){35}} \put(50,4){\vector(-1,0){35}}
\curve(50,6, 50,2) \put(47,0){$\Sigma$} \put(78,5){$v_1$}
\put(18,5){$v_2$} \put(65,-2){Figure 1}
\end{picture}
\end{figure}
The mirror manifold $Y$, as a symplectic manifold, is the cylinder
$Y=\mathbb{R}\times\sqrt{-1}S^1$. Any Lagrangian section $L$ of
$\mu_Y:Y\rightarrow\mathbb{R}$ lifts to the universal cover
$T^*N_\mathbb{R}=\mathbb{R}^2$. A lift $\tilde{L}$ of $L$ is given
by a graph
$$\tilde{L}=\{(\xi,y(\xi)):\xi\in\mathbb{R}\}\subset\mathbb{R}^2,$$
where $y(\xi)=g'(\xi)$ is the derivative of a function $g=g(\xi)\in
C^2(\mathbb{R})$. Given $(a,b)\in\mathbb{Z}^2$, the conditions in
Definition~\ref{def3.1} reduces to the following two equalities of
limits:
\begin{eqnarray*}
\lim_{\xi\rightarrow-\infty}2e^{-2\xi}(y(\xi)+a) & = &
\lim_{\xi\rightarrow-\infty}e^{-2\xi}y'(\xi),\\
\lim_{\xi\rightarrow\infty}2e^{2\xi}(b-y(\xi)) & = &
\lim_{\xi\rightarrow\infty}e^{2\xi}y'(\xi).
\end{eqnarray*}
This implies that, geometrically, we have $y(\xi)\rightarrow-a$ as
$\xi\rightarrow-\infty$ and $y(\xi)\rightarrow b$ as
$\xi\rightarrow\infty$, and the slope of the graph goes to zero as
$\xi\rightarrow\pm\infty$; there are no restrictions on the graph
for finite values of $\xi$. The equalities of limits place further
restrictions on the growth rates of $y(\xi)$ and its derivative as
$\xi$ tends to $\pm\infty$.\\

\noindent\textbf{Example 2.} Consider the case when
$X_\Sigma=\mathbb{C}P^2$. The fan $\Sigma$ in
$N_\mathbb{R}=\mathbb{R}^2$ is generated by
$v_1=(1,0),v_2=(0,1),v_3=(-1,-1)$ (see Figure 2 below).
\begin{figure}[ht]
\setlength{\unitlength}{1mm}
\begin{picture}(100,30)
\put(50,12){\vector(1,0){20}} \put(50,12){\vector(0,1){20}}
\put(50,12){\vector(-1,-1){15}} \put(70,13){$v_1$}
\put(51,31){$v_2$} \put(34,0){$v_3$} \put(57,20){$\sigma_1$}
\put(40,15){$\sigma_2$} \put(52,4){$\sigma_3$} \put(70,-2){Figure 2}
\end{picture}
\end{figure}
The mirror manifold $Y$ is given by
$Y=\mathbb{R}^2\times\sqrt{-1}T^2$ equipped with the standard
symplectic structure. Any Lagrangian section $L$ of
$\mu_Y:Y\rightarrow\mathbb{R}^2$ can be lifted to a graph
$$\tilde{L}=\{(\xi_1,\xi_2,y_1(\xi_1,\xi_2),y_2(\xi_1,\xi_2)):(\xi_1,\xi_2)\in\mathbb{R}^2\}$$
in the universal cover $T^*N_\mathbb{R}=\mathbb{R}^4$, where
$y_1(\xi_1,\xi_2)=\frac{\partial
g}{\partial\xi_1},y_2(\xi_1,\xi_2)=\frac{\partial g}{\partial\xi_2}$
are the partial derivatives of a function $g=g(\xi_1,\xi_2)\in
C^2(\mathbb{R}^2)$. Let $(a_1,a_2,a_3)\in\mathbb{Z}^3$. Consider the
maximal cone $\sigma_1$. Then the conditions in
Definition~\ref{def3.1} can be restated as
\begin{eqnarray*}
\lim_{\xi_1\rightarrow-\infty}2e^{-2\xi_1}(y_1(\xi_1,\xi_2)+a_1) & =
& \lim_{\xi_1\rightarrow-\infty}e^{-2\xi_1}y_{1,1}(\xi_1,\xi_2),\\
\lim_{\xi_2\rightarrow-\infty}2e^{-2\xi_2}(y_2(\xi_1,\xi_2)+a_2) & =
& \lim_{\xi_2\rightarrow-\infty}e^{-2\xi_2}y_{2,2}(\xi_1,\xi_2),\\
\lim_{\xi_1\rightarrow-\infty}e^{-\xi_1-\xi_2}y_{1,2}(\xi_1,\xi_2) &
= &
\lim_{\xi_1\rightarrow-\infty}e^{-\xi_1-\xi_2}y_{1,2}(\xi_1,\xi_2)=0,
\end{eqnarray*}
where we denote by $y_{i,j}$ the partial derivative $\frac{\partial
y_i}{\partial\xi_j}$. In particular, we must have
$y_1\rightarrow-a_1$ as $\xi_1\rightarrow-\infty$,
$y_2\rightarrow-a_2$ as $\xi_2\rightarrow-\infty$, and various
partial derivatives of $y_1,y_2$ go to zero as $t_1,t_2$ tends to
$-\infty$.

For another maximal cone, say, $\sigma_2$, the conditions can
similarly be rewritten as
\begin{eqnarray*}
\lim_{\xi_2\rightarrow-\infty}2e^{-2\xi_2}(y_2(\xi_1,\xi_1+\xi_2)+a_2)
& = &
\lim_{\xi_2\rightarrow-\infty}e^{-2\xi_2}y_{2,2}(\xi_1,\xi_1+\xi_2),\\
\lim_{\xi_1\rightarrow\infty}2e^{2\xi_1}((-y_1-y_2)(\xi_1,\xi_1+\xi_2)+a_3)
& = &
\lim_{\xi_1\rightarrow\infty}e^{2\xi_1}(y_{1,1}+2y_{1,2}+y_{2,2})(\xi_1,\xi_1+\xi_2),\\
\lim_{\xi_1\rightarrow\infty}e^{\xi_1-\xi_2}(y_{1,2}+y_{2,2})(\xi_1,\xi_1+\xi_2)
& = &
\lim_{\xi_2\rightarrow-\infty}e^{\xi_1-\xi_2}(y_{1,2}+y_{2,2})(\xi_1,\xi_1+\xi_2)=0.
\end{eqnarray*}
Geometrically, this means that we should also have
$-y_1-y_2\rightarrow-a_3$ as $\xi$ goes to $-\infty$ in the
$(-1,-1)$ direction, and various combinations of the partial
derivatives of $y_1,y_2$ go to zero as $\xi$ tends to $-\infty$ in
the $(0,1)$ and $(-1,-1)$ directions. Again, the equalities of
limits indicate the growth rates of $y_1,y_2$ and combinations of
their partial derivatives as $\xi$ tends to $-\infty$ in the
$(1,0),(0,1),(-1,-1)$ directions.\\

In general, given a lift $\tilde{L}=\{(\xi,y(\xi)):\xi\in
N_\mathbb{R}\}$ of a Lagrangian section $L$ of from $Y$ to
$T^*N_\mathbb{R}$, the conditions in Definition~\ref{def3.1} specify
the values and growth rates of the functions
$y_1(\xi),\ldots,y_n(\xi)$ and combinations of their partial
derivatives as $\xi$ tends to $-\infty$ in the directions of
$v_1,\ldots,v_d$. In particular, for $i=1,\ldots,d$, $\langle
y(\xi),v_i\rangle$ goes to $-a_i$ as $\xi$ tends to $-\infty$ in the
direction of $v_i$.

We may also regard $L$ and any lift $\tilde{L}$ of $L$ as Lagrangian
sections over $P$, the interior of the polytope $\bar{P}$, so that
we can write $\tilde{L}=\{(x,y(x)):x\in P\}$. Then the conditions in
Definition~\ref{def3.1} can be viewed as boundary conditions for the
functions $y_1(x),\ldots,y_n(x)$ and combinations of their partial
derivatives over the boundary $\partial\bar{P}$. For example, if
$\tilde{L}$ satisfies $(\ast_a)$, where
$a=(a_1,\ldots,a_d)\in\mathbb{Z}^d$, then the function $\langle
y(x),v_k\rangle$ tends to $-a_k$ as $x$ approaches the facet of
$\bar{P}$ with normal vector $v_k$.
\begin{nb}
Our Lagrangian sections are closely related to the
\textit{tropical Lagrangian sections} defined and used by Abouzaid
in his proof \cite{Abouzaid05}, \cite{Abouzaid06} of the Homological
Mirror Symmetry Conjecture for toric varieties. This relation is
similar to the one explained in Appendix C of
Fang-Liu-Treumann-Zaslow \cite{FLTZ08}.\footnote{Indeed, we believe
that the boundary conditions for Lagrangian sections used in
\cite{FLTZ08} are equivalent to those we use here.} Let us describe
the relation briefly as follows. In \cite{Abouzaid05},
\cite{Abouzaid06}, Abouzaid considered the family of superpotentials
$$W_t=\sum_{i=1}^d c_it^{-\lambda_i}z^{v_i},$$
and the smooth hypersurfaces $M_t=W_t^{-1}(0)$ in
$TM_\mathbb{R}/M=(\mathbb{C}^*)^n$. The \textit{amoeba} of $M_t$ is
the image under the logarithm map, i.e.
$\mathcal{A}_t=\textrm{Log}(M_t)\subset M_\mathbb{R}$, and the
tropical amoeba is the limit
$\Pi=\lim_{t\rightarrow\infty}(\mathcal{A}_t/\log t)\subset
M_\mathbb{R}$. Abouzaid showed that there is a distinguished
connected component $Q$ of $M_\mathbb{R}\setminus\Pi$ which is a
copy of the moment polytope $\bar{P}$ of
$(X_\Sigma,\omega_{X_\Sigma})$. Abouazid then defined his tropical
Lagrangian sections to be the Lagrangian sections over $Q$ with
boundary in $M_\infty=\lim_{t\rightarrow\infty}M_t$. Now, given a
Lagrangian section $L$ in $Y$ satisfying $(\ast_{[a]})$ for some
$[a]\in A(W)$, we may regard $L$ as a Lagrangian section over $P$
(by writing $Y$ as $P\times\sqrt{-1}T_M$) and hence over the
interior of $Q\subset M_\mathbb{R}\setminus\Pi$. Then $L$ is in the
equivalence class of Abouzaid's Lagrangian section associated to the
line bundle $\mathcal{L}_{[a]}$.
\end{nb}
We now return to the general discussion of the set
$\mathbbm{L}(Y,W)$ of Lagrangian sections.
\begin{prop}
Any two Lagrangian sections $L_1,L_2\in\mathbbm{L}(Y,W)$ satisfying
the same growth condition $(\ast_{[a]})$ can be deformed to each
other through Hamiltonian isotopies which preserve $(\ast_{[a]})$.
\end{prop}
\begin{proof}
Choose lifts $\tilde{L}_1,\tilde{L}_2$ of $L_1,L_2$ respectively,
such that they satisfy the same growth condition $(\ast_a)$, for
some representative $a\in\mathbb{Z}^d$ of $[a]$. Let $g_1, g_2$ be
the potentials of $\tilde{L}_1,\tilde{L}_2$ respectively. Regard
$H:=g_1-g_2$ as a $T_M$-invariant function on $Y$. Then the
Hamiltonian flow $\rho_t:Y\rightarrow Y$ associated to $H$ moves
$L_1$ to $L_2$ at time $t=1$, and $\rho_t(L_1)$ satisfies
$(\ast_{[a]})$ for all $t$ because $H$, as a function on
$N_\mathbb{R}$, satisfies $(\ast_0)$.
\end{proof}
In view of this proposition, we define two Lagrangian sections
$L_1,L_2\in\mathbbm{L}(Y,W)$ to be equivalent, denoted $L_1\sim
L_2$, if they satisfy the same growth condition $(\ast_{[a]})$; and
we denote the equivalence class to which $L\in\mathbbm{L}(Y,W)$
belongs by $[L]$.\\

Now rewrite $Y$ as $Y=P\times\sqrt{-1}T_M=TP/M$ and use the
coordinates $x_j$'s and $y_j$'s to express a lift
$\tilde{L}=\{(\xi,dg(\xi)):\xi\in N_\mathbb{R}\}$ of the Lagrangian
section $L$ as the graph of the gradient of the function $\Psi^*g$,
with respect to the metric $\sum_{j,k=1}^n\psi_{jk}dx_j\otimes dx_k$
on $P$. In other words, we have $\tilde{L}=\{(x,y(x)):x\in P,
y(x)=\nabla(\Psi^*g)(x)\}$, or in coordinates,
$$y_j(x)=\sum_{k=1}^n\psi^{jk}\frac{\partial(\Psi^*g)}{\partial x_k}.$$
For any Lagrangian section $L$ of $\nu_Y:Y\rightarrow P$, define the
\textit{normalized slope} of $L$ by
$$\lambda(L)=\frac{1}{\textrm{Vol}(P)}\int_P\sum_{j=1}^n
\frac{\partial y_j(x)}{\partial x_j}dx_1\wedge\ldots\wedge dx_n,$$
where $\tilde{L}=\{(x,y(x)):x\in P\}\subset T^*N_\mathbb{R}$ is any
lift of $L$ to $T^*N_\mathbb{R}$. $\lambda(L)$ is clearly
independent of the choice of the lift $\tilde{L}$.
\begin{prop}
If $L_1\sim L_2$, then $\lambda(L_1)=\lambda(L_2)$. Hence $\lambda$
is an invariant on the set of equivalence classes
$\mathbbm{L}(Y,W)/\sim$.
\end{prop}
\begin{proof}
As in the proof of the above proposition, we choose lifts
$\tilde{L}_1,\tilde{L}_2$ of $L_1,L_2$ respectively such that they
satisfy the same growth condition $(\ast_a)$, for some
$a\in\mathbb{Z}^d$ representing $[a]$. Let $g_1, g_2$ be the
potentials of $\tilde{L}_1,\tilde{L}_2$ respectively, and let
$H:=g_1-g_2$. Set
$y_j(x)=\sum_{k=1}^n\psi^{jk}\frac{\partial(\Psi^*H)}{\partial
x_k}$. Then, for $j=1,\ldots,n$, we have
\begin{eqnarray*}
\int_P\sum_{j=1}^n\frac{\partial y_j}{\partial
x_j}dx_1\wedge\ldots\wedge dx_n & = &
\int_Pd(\sum_{j=1}^n(-1)^{j-1}y_jdx_1\wedge\ldots\wedge\widehat{dx_j}\wedge\ldots\wedge
dx_n)\\
& = & \int_{\partial\bar{P}}
\sum_{j=1}^n(-1)^{j-1}y_jdx_1\wedge\ldots\wedge\widehat{dx_j}\wedge\ldots\wedge
dx_n,
\end{eqnarray*}
by Stokes theorem. Consider a facet $F_k=\{x\in\bar{P}:l_k(x)=0\}$
of $\bar{P}$. Without loss of generality, suppose that $v_k^n\neq0$.
Then use $x_1,\ldots,x_n$ as the coordinates on $F_k$, so that
$x_n=-\frac{\lambda_k}{v_k^n}-\sum_{p=1}^{n-1}\frac{v_k^p}{v_k^n}x_p$.
We have
$$\sum_{j=1}^n(-1)^{j-1}y_jdx_1\wedge\ldots\wedge\widehat{dx_j}\wedge\ldots\wedge
dx_n=\frac{(-1)^{n-1}}{v_k^n}\langle y(x),v_k\rangle
dx_1\wedge\ldots\wedge dx_{n-1}.$$ Now, since $H$ satisfies
$(\ast_0)$, $\langle y(x),v_k\rangle=0$ for $x\in F_k$. Hence
$\int_P\sum_{j=1}^n\frac{\partial y_j}{\partial
x_j}dx_1\wedge\ldots\wedge dx_n=0$, and we have
$\lambda(L_1)=\lambda(L_2)$.
\end{proof}
\begin{defn}
A Lagrangian section $L\in\mathbbm{L}(Y,W)$ is said to be harmonic
if the following Laplace-type equation is satisfied
\begin{equation}\label{harmonic}
\sum_{j=1}^n \frac{\partial y_j(x)}{\partial x_j}=\lambda(L),
\end{equation}
for some lift $\tilde{L}=\{(x,y(x)):x\in P\}\subset T^*N_\mathbb{R}$
of $L$.
\end{defn}
\noindent The equation (\ref{harmonic}) is equivalent to the
following equation
$$\sum_{j,k=1}^n\psi^{jk}\Bigg(\frac{\partial^2(\Psi^*g)}{\partial x_j\partial x_k}
-\sum_{p,q=1}^n\psi^{pq}\psi_{pjk}\frac{\partial(\Psi^*g)}{\partial
x_q}\Bigg)=\lambda(L)$$ on $P$, where $\psi_{pjk}$ denotes
$\frac{\partial^3\psi}{\partial x_p\partial x_j\partial x_k}$. If we
regard $L=\{(\xi,dg(\xi)):\xi\in N_\mathbb{R}\}$ as a section of
$\mu_Y:Y\rightarrow N_\mathbb{R}$, then $L$ is harmonic if and only
of $g$ is a solution to the equation
$$\sum_{j,k=1}^n\phi^{jk}\frac{\partial^2g}{\partial\xi_j\partial\xi_k}=\lambda(L)$$
on $N_\mathbb{R}$. In the next section, we will see that in each
equivalence class $[L]\in\mathbbm{L}(Y,W)/\sim$ of Lagrangian
sections, there exists a unique harmonic representative. This is
mirror to the existence of a unique Hermitian-Einstein metric on
each holomorphic line bundle over $X_\Sigma$, and $\lambda(L)$ is
the mirror analogue of the (normalized) slope of a line bundle.

On the other hand, we may also choose special Lagrangian sections as
representatives. According to the definition of
Auroux~\cite{Auroux07}, a Lagrangian submanifold $L\subset Y$ is
special with phase $\theta\in\mathbb{R}$ if
$\textrm{Im}(e^{\sqrt{-1}\theta}\Omega_Y)|_L=0$. In terms of the
$x_j,y_j$ coordinates,
\begin{eqnarray*}
\Omega_Y|_L &=&
\bigwedge_{j=1}^n\Bigg(-dx_j+\sqrt{-1}dy_j(x)\Bigg)\\
&=&
\bigwedge_{j=1}^n\Bigg(\sum_{k=1}^n\Bigg(-\delta_{jk}+\sqrt{-1}\frac{\partial
y_j(x)}{\partial x_k}\Bigg)dx_k\Bigg)\\
&=& \textrm{det}\Bigg(-I_n+\sqrt{-1}\Big(\frac{\partial
y_j(x)}{\partial x_k}\Big)_{j,k=1}^n\Bigg)dx_1\wedge\ldots\wedge
dx_n,
\end{eqnarray*}
where $I_n$ denotes the $n\times n$ identity matrix. So
$L=\{(x,y(x):x\in P\}$ is special Lagrangian with phase
$\theta\in\mathbb{R}$ if and only if the following equation is
satisfied
\begin{equation}\label{slag}
\textrm{Im}\Bigg(e^{\sqrt{-1}\theta}\textrm{det}\Bigg(I_n-\sqrt{-1}\Big(\frac{\partial
y_j(x)}{\partial x_k}\Big)_{j,k=1}^n\Bigg)\Bigg)=0.
\end{equation}
Equivalently, this means $\Psi^*g$ satisfies the equation
$$\textrm{Im}\Bigg(e^{\sqrt{-1}\theta}\textrm{det}\Bigg(I_n-\sqrt{-1}
\Bigg[\sum_{l=1}^n\psi^{jl}\Bigg(\frac{\partial^2(\Psi^*g)}{\partial
x_l\partial x_k} -\sum_{p,q=1}^n\psi^{pq}\psi_{plk}\frac{\partial
(\Psi^*g)}{\partial x_q}\Bigg)\Bigg]_{j,k=1}^n\Bigg)\Bigg)=0.$$ or,
in the $\xi_j,y_j$ coordinates, $g$ satisfies the equation
$$\textrm{Im}\Bigg(e^{\sqrt{-1}\theta}\textrm{det}\Bigg(I_n-\sqrt{-1}
\Bigg(\sum_{l=1}^n\phi^{kl}\frac{\partial^2g}{\partial\xi_j\partial\xi_l}
\Bigg)_{j,k=1}^n\Bigg)\Bigg)=0.$$ Our harmonic Lagrangians are
closely related to special Lagrangians, at least in the large radius
limit: If we rescale the fiber coordinates by replacing $y_j$ by
$\epsilon y_j$, then, for small $\epsilon$, the leading term of
equation (\ref{slag}) will give
$$\sum_{j=1}^n \frac{\partial y_j(x)}{\partial x_j}=\frac{1}{\epsilon}\tan\theta,$$
which is nothing but equation (\ref{harmonic}) if we choose $\theta$
such that $\tan\theta=\epsilon\lambda(L)$.

\section{The SYZ mirror transformation as a geometric
correspondence}\label{sec4}

In this section, we first recall the definition of the SYZ mirror
transformation. Then we proceed to prove our main result.

For $[a]\in H^2(X_\Sigma,\mathbb{Z})$, let $\mathcal{L}_{[a]}$ be
the corresponding holomorphic line bundle over $X_\Sigma$. Choose a
$T_N$-equivariant meromorphic section $s$ of $\mathcal{L}_{[a]}$.
Then $\textrm{div}(s)=\sum_{i=1}^da_iD_i$, for some integers
$a_1,\ldots,a_d\in\mathbb{Z}$ such that
$(a_1,\ldots,a_d)\in\mathbb{Z}^d$ gives a representative of the
class $a$. Note that $s$ is holomorphic and nowhere vanishing over
$U\subset X_\Sigma$, so it is a holomorphic frame of
$\mathcal{L}_{[a]}|_U$.

Let $h$ be a $T_N$-invariant hermitian metric of class $C^2$ on
$\mathcal{L}_{[a]}$. The Chern connection $\nabla_h$ is given by
$\nabla_h=d+\partial\log h(s,s)$ over $U$. If we define a function
$g_h:N_\mathbb{R}\rightarrow\mathbb{R}$ by setting
$$g_h(\xi)=-\frac{1}{2}\log
h(s(e^{\xi+\sqrt{-1}u}),s(e^{\xi+\sqrt{-1}u})),$$ then the
restriction of $\nabla_h$ to a fiber $F_\xi:=\nu_U^{-1}(\xi)\cong
T_N$ gives a flat $U(1)$-connection
$$d+\frac{\sqrt{-1}}{2}\sum_{j=1}^n\frac{\partial\log
h(s,s)}{\partial\xi_j}du_j=d-\sqrt{-1}\sum_{j=1}^n\frac{\partial
g_h}{\partial\xi_j}du_j$$ on the trivial line bundle
$\underline{\mathbb{C}}$ over $T_N$. Recall that the dual torus
$T_M=(T_N)^*$ can be interpreted as the space of flat
$U(1)$-connections on the trivial line bundle
$\underline{\mathbb{C}}$ over $T_N$ modulo gauge
equivalence.\footnote{This is in fact the starting point of the SYZ
conjecture \cite{SYZ96}} In our situation, the connection
$d-\sqrt{-1}\sum_{j=1}^n\frac{\partial g_h}{\partial\xi_j}du_j$
corresponds to the point $(\frac{\partial
g_h}{\partial\xi_1},\ldots,\frac{\partial g_h}{\partial\xi_n})\in
T_M$. Hence, the hermitian metric $h$, or the Chern connection
$\nabla_h$, determines a section
$$\tilde{L}_h=\{(\xi,dg_h(\xi))=(\xi_1,\ldots,\xi_n,\frac{\partial
g_h}{\partial\xi_1},\ldots,\frac{\partial
g_h}{\partial\xi_n}):\xi\in N_\mathbb{R}\}$$ of
$T^*N_\mathbb{R}=N_\mathbb{R}\times\sqrt{-1}M_\mathbb{R}\rightarrow
N_\mathbb{R}$, which is Lagrangian since $\nabla_h$ is holomorphic
(see \cite{LYZ00}). $\tilde{L}_h$ descends to give a Lagrangian
section $L_h$ of $\mu_Y:Y\rightarrow N_\mathbb{R}$.

If $s'$ is another $T_N$-equivariant meromorphic section of
$\mathcal{L}_a$, then $s'=cw^u\cdot s$, for some constant
$c\in\mathbb{C}^*$ and $u\in M$, where $w^u$ is the monomial
$w_1^{u^1}\ldots w_n^{u^n}$. Since
$h(s'(w),s'(w))=|cw^u|^2h(s(w),s(w))=|c|^2e^{2\langle
u,\xi\rangle}h(s(w),s(w))$, we have $g_h'(\xi)=-\log |c|-\langle
u,\xi\rangle+g_h(\xi)$, where $g_h':=-\frac{1}{2}\log h(s',s')$. So
$dg_h'(\xi)=dg_h(\xi)-u$. This gives a different Lagrangian section
$\tilde{L}'_h=\{(\xi,dg_h'(\xi)):\xi\in N_\mathbb{R}\}$ in
$T^*N_\mathbb{R}$, but it descends to the same Lagrangian section
$L_h$ in $Y$.

Thus we have a well-defined transformation
\begin{equation*}
\mathcal{F}:h\mapsto L_h
\end{equation*}
from the set of $T_N$-invariant hermitian metrics on holomorphic
line bundles over $X_\Sigma$ to the set of Lagrangian sections of
$\mu_Y:Y\rightarrow N_\mathbb{R}$. This is called the \textit{SYZ
mirror transformation}. This is (fiberwise) a real version of the
Fourier-Mukai transform in algebraic geometry. We can invert the
construction and define the inverse SYZ mirror transformation
$\mathcal{F}^{-1}$, which produces, from a Lagrangian section $L$ of
$\mu_Y:Y\rightarrow N_\mathbb{R}$, a $T_N$-invariant hermitian
metric $h_L:=\mathcal{F}^{-1}(L)$ on a holomorphic line bundle
\textit{over $U$}. However, $h_L$ may not be extended to a hermitian
metric on a holomorphic line bundle over $X_\Sigma$. The question we
raised in the introduction is to characterize the set of Lagrangian
sections $L$ for which $h_L$ can be extended over $X_\Sigma$. Our
main result says that this set is precisely $\mathbbm{L}(Y,W)$,
which we introduced in the last section.
\begin{thm}\label{main}
The image of the SYZ mirror transformation $\mathcal{F}$ is
$\mathbbm{L}(Y,W)$, i.e. for a Lagrangian section $L$ of
$\mu_Y:Y\rightarrow N_\mathbb{R}$, there exists a $T_N$-invariant
hermitian metric $h$ on a holomorphic line bundle over $X_\Sigma$
such that $L=L_h=\mathcal{F}(h)$ if and only if $L$ satisfies the
growth condition $(\ast_{[a]})$ for some $[a]\in A(W)$.
\end{thm}
Before we prove the theorem, we need a couple of lemmas. Let $[a]$
be an element in $A(W)=H^2(X_\Sigma,\mathbb{Z})$ and
$\mathcal{L}_{[a]}$ the corresponding holomorphic line bundle over
$X_\Sigma$. We first consider a particular $T_N$-invariant hermitian
metric $h_0$ on $\mathcal{L}_{[a]}$ defined as follows. Choose a
representative $(a_1,\ldots,a_n)\in\mathbb{Z}^d$ of $a$, and fix a
$T_N$-equivariant meromorphic section $s$ of $\mathcal{L}_{[a]}$
such that $\textrm{div}(s)=D_a=\sum_{i=1}^da_iD_i$, so that we can
canonically identify $\mathcal{L}_{[a]}$ with the toric divisor line
bundle $\mathcal{O}(D_a)$. Recall that the moment map
$\mu_U:U\rightarrow P$ is given by
$$\mu_U(w)=d\phi(\log|w_1|,\ldots,\log|w_n|)=
\frac{\sum_{u\in\bar{P}\cap M}c_u|w^u|^2\cdot
u}{\sum_{u\in\bar{P}\cap M}c_u|w^u|^2},$$ for $w=(w_1,\ldots,w_n)\in
U=(\mathbb{C}^*)^n$. For $i=1,\ldots,d$, let
$l_i:M_\mathbb{R}\rightarrow\mathbb{R}$ be the function defined by
$l_i(x)=\langle x,v_i\rangle+\lambda_i$. In \cite{Guillemin94b},
Guillemin showed that there is a $T_N$-invariant hermitian metric
$h_0$ on $\mathcal{L}_{[a]}$ such that
$$h_0(s,s)=\prod_{i=1}^d(l_i\circ\mu_U)^{a_i}.$$
\begin{lem}\label{lem1}
$\tilde{L}_{h_0}$ satisfies the growth condition $(\ast_a)$
\end{lem}
The proof of this lemma, which is a straightforward but lengthy
calculation, will be given in the appendix.

To describe the other lemma we require, consider the diagonal
$T^n$-action on $\mathbb{C}^n$. If
$F:\mathbb{C}^n\rightarrow\mathbb{R}$ is a $T^n$-invariant function,
then we can define a function $f:\mathbb{R}^n\rightarrow\mathbb{R}$,
by
$f(\xi_1,\ldots,\xi_n)=F(e^{\xi_1+\sqrt{-1}u_1},\ldots,e^{\xi_n+\sqrt{-1}u_n})$,
where $w_j=e^{\xi_j+\sqrt{-1}u_j}$, $j=1,\ldots,n$, are the complex
coordinates on $\mathbb{C}^n$. But not all functions on
$\mathbb{R}^n$ come from this way.
\begin{lem}\label{lem2}
Given a function $f\in C^2(\mathbb{R}^n)$. Define
$F:(\mathbb{C}^*)^n\rightarrow\mathbb{R}$ by
$F(w_1,\ldots,w_n)=f(\log|w_1|,\ldots,\log|w_n|)$. Then $F$ can be
extended to a $T^n$-invariant $C^2$ function on $\mathbb{C}^n$ if
and only if the following three conditions are satisfied
\begin{enumerate}
\item[1.] For $j=1,\ldots,n$,
$e^{-2\xi_j}\frac{\partial^2f}{\partial\xi_j^2}$ and
$2e^{-2\xi_j}\frac{\partial f}{\partial\xi_j}$ go to the same limit
as $\xi_j\rightarrow-\infty$.
\item[2.] For any $j,k,l\in\{1,\ldots,n\}$, the limit of
$\frac{\partial^2f}{\partial\xi_j\partial\xi_k}$ exists as
$\xi_l\rightarrow-\infty$.
\item[3.] For any distinct $j,k\in\{1,\ldots,n\}$,
$e^{-\xi_j-\xi_k}\frac{\partial^2f}{\partial\xi_j\partial\xi_k}$
goes to zero as $\xi_j\rightarrow-\infty$ or
$\xi_k\rightarrow-\infty$.
\end{enumerate}
\end{lem}
\begin{proof}
Write $e^{\xi_j+\sqrt{-1}u_j}=w_j=x_j+\sqrt{-1}y_j$. Then, by the
chain rule, we have, for $j=1,\ldots,n$,
\begin{eqnarray*}
\frac{\partial F}{\partial x_j} & = & e^{-\xi_j}\cos
u_j\frac{\partial f}{\partial\xi_j},\ \frac{\partial F}{\partial
y_j}=e^{-\xi_j}\sin u_j\frac{\partial f}{\partial\xi_j},\\
\frac{\partial^2F}{\partial x_j^2} & = &
e^{-2\xi_j}\cos^2u_j(\frac{\partial^2f}{\partial\xi_j^2}-2\frac{\partial
f}{\partial\xi_j})+e^{-2\xi_j}\frac{\partial f}{\partial\xi_j},\\
\frac{\partial^2F}{\partial x_j\partial y_j} & = & e^{-2\xi_j}\cos
u_j\sin u_j(\frac{\partial^2f}{\partial\xi_j^2}-2\frac{\partial
f}{\partial\xi_j}),\\
\frac{\partial^2F}{\partial y_j^2} & = &
e^{-2\xi_j}\sin^2u_j(\frac{\partial^2f}{\partial\xi_j^2}-2\frac{\partial
f}{\partial\xi_j})+e^{-2\xi_j}\frac{\partial f}{\partial\xi_j},
\end{eqnarray*}
and, for $j\neq k$,
\begin{eqnarray*}
\frac{\partial^2F}{\partial x_j\partial x_k} & = &
e^{-\xi_j-\xi_k}\cos u_j\cos
u_k\frac{\partial^2f}{\partial\xi_j\partial\xi_k},\\
\frac{\partial^2F}{\partial x_j\partial y_k} & = &
e^{-\xi_j-\xi_k}\cos u_j\sin
u_k\frac{\partial^2f}{\partial\xi_j\partial\xi_k},\\
\frac{\partial^2F}{\partial y_j\partial y_k} & = &
e^{-\xi_j-\xi_k}\sin u_j\sin
u_k\frac{\partial^2f}{\partial\xi_j\partial\xi_k}.
\end{eqnarray*}
It is then not hard to see that the conditions (1)-(3) are necessary
and sufficient conditions for extending $F$ to $\mathbb{C}^n$.
\end{proof}
\begin{proof}[Proof of Theorem~\ref{main}]
Let $h$ be any other $T_N$-invariant $C^2$ hermitian metric of
$\mathcal{L}_{[a]}$. Then there is a function $F\in C^2(X_\Sigma)$
such that $h=e^{-2F}h_0$. Restrict $F$ to $U\subset X_\Sigma$, and
define $f:N_\mathbb{R}\rightarrow\mathbb{R}$ by
$f(\xi_1,\ldots,\xi_n)=F(e^{\xi_1+\sqrt{-1}u_1},\ldots,e^{\xi_n+\sqrt{-1}u_n})$.
Let $\sigma\in\Sigma$ be an $n$-dimensional cone, and
$U_\sigma=\textrm{Spec }\mathbb{C}[\check{\sigma}\cap M]$ the
corresponding affine toric variety. $X_\Sigma$ is covered by these
$U_\sigma$'s, and since $X_\Sigma$ is nonsingular,
$U_\sigma\cong\mathbb{C}^n$. Without loss of generality, suppose
that the generators of $\sigma$ are $v_1,\ldots,v_n\in N$. They give
a $\mathbb{Z}$-basis of $N$. Let
$\tilde{w}_1=e^{\tilde{\xi}_1+\sqrt{-1}\tilde{u}_1},\ldots,\tilde{w}_n
=e^{\tilde{\xi}_n+\sqrt{-1}\tilde{u}_n}$ be the corresponding
(inhomogeneous) complex coordinates on $U_\sigma$. This gives
coordinates $\tilde{\xi}_1,\ldots,\tilde{\xi}_n$ on $N_\mathbb{R}$,
and the transformation from these coordinates to the original
coordinates $\xi_1,\ldots,\xi_n$ is given by
$$\xi=(\xi_1,\ldots,\xi_n)=v_1\tilde{\xi}_1+\ldots+v_n\tilde{\xi}_n.$$
Apply the chain rule, we get
\begin{eqnarray*}
\frac{\partial f}{\partial\tilde{\xi}_i} & = &
\sum_{k=1}\frac{\partial
f}{\partial\xi_k}\frac{\partial\xi_k}{\partial\tilde{\xi}_i}=\sum_{k=1}v_i^k\frac{\partial
f}{\partial\xi_k}=\langle df,v_i\rangle,\\
\frac{\partial^2f}{\partial\tilde{\xi}_j\partial\tilde{\xi}_k} & =
& \sum_{p,q=1}^nv_j^pv_k^q
\frac{\partial^2f}{\partial\xi_p\partial\xi_q}=v_j^T\textrm{Hess}(f)v_k.
\end{eqnarray*}
Hence, by Lemma~\ref{lem2}, we conclude that the function $f$
satisfies the growth condition $(\ast_0)$. Now, by Lemma~\ref{lem1},
$\tilde{L}_{h_0}$ satisfies the growth condition $(\ast_a)$. Since
$g_h=g_{h_0}+f$, we see that $\tilde{L}_h$ also satisfies
$(\ast_a)$.

Conversely, let $L$ be a Lagrangian section in $Y$ satisfying
$(\ast_{[a]})$. Choose a lift $\tilde{L}=\{(\xi,dg(\xi)):\xi\in
N_\mathbb{R}\}\subset T^*N_\mathbb{R}$ of $L$ which satisfies the
same growth condition $(\ast_a)$ as $\tilde{L}_{h_0}$. Then the
$C^2$ function $f:=g-g_{h_0}:N_\mathbb{R}\rightarrow\mathbb{R}$
satisfies the growth condition $(\ast_0)$. By the above argument and
Lemma~\ref{lem2}, $f$ extends to a function $F\in C^2(X_\Sigma)$. So
$h:=e^{-2F}h_0$ defines a $T_N$-invariant hermitian metric on
$\mathcal{L}_{[a]}$. This completes the proof of the theorem.
\end{proof}
Theorem~\ref{main} establishes a bijective correspondence between
$T_N$-invariant hermitian metrics on the holomorphic line bundle
$\mathcal{L}_{[a]}$ over $X_\Sigma$ and Lagrangian sections of
$(Y,W)$ satisfying the growth condition $(\ast_{[a]})$, for any
$[a]\in H^2(X_\Sigma,\mathbb{Z})=A(W)$. In addition, by our
definition in Section~\ref{sec3}, two Lagrangian sections
$L_1,L_2\in\mathbbm{L}(Y,W)$ are equivalent, denoted $L_1\sim L_2$,
if and only if they satisfies the same growth condition. Hence, an
immediate consequence of our main result is the following
\begin{cor}
The SYZ mirror transformation $\mathcal{F}$ induces a bijective map
$$\mathcal{F}:\textrm{Pic}(X_\Sigma)\overset{\cong}{\rightarrow}(\mathbbm{L}(Y,W)/\sim).$$
\end{cor}
Recall that a hermitian metric $h$ on the line bundle
$\mathcal{L}_{[a]}$ is Hermitian-Einstein, with respect to the
K\"{a}hler metric $\omega_{X_\Sigma}$ on $X_\Sigma$, if and only if
the following equation is satisfied
$$\sqrt{-1}F_h\wedge\omega_{X_\Sigma}^{n-1}=\frac{\lambda(\mathcal{L}_{[a]})}{n}\cdot\omega_{X_\Sigma}^n,$$
where $F_h$ is the curvature of the Chern connection $\nabla_h$, and
$\lambda(\mathcal{L}_{[a]})$ is the \textit{normalized slope} of
$\mathcal{L}_{[a]}$ defined by
$$\lambda(\mathcal{L}_{[a]}):=\frac{n\cdot\int_{X_\Sigma}\sqrt{-1}F_h\wedge\omega_{X_\Sigma}^{n-1}}
{\int_{X_\Sigma}\omega_{X_\Sigma}^n}=\frac{2\pi
n\mu(\mathcal{L}_{[a]})}{\int_{X_\Sigma}\omega_{X_\Sigma}^n}.$$ Now
let $y_j=\frac{\partial
g_h}{\partial\xi_j}=\sum_{k=1}^n\psi^{jk}\frac{\partial\Psi^*g_h}{\partial
x_k}$, then, restricting to $U\subset X_\Sigma$, we have
$$\sqrt{-1}F_h=\bar{\partial}\partial\log h=\sum_{j=1}^ndy_j\wedge
du_j.$$ Hence,
\begin{eqnarray*}
\sqrt{-1}F_h\wedge\omega_{X_\Sigma}^{n-1} & = &
(\sum_{j=1}^ndy_j\wedge du_j)\wedge(\sum_{j=1}^ndx_j\wedge
du_j)^{n-1}\\
& = & (n-1)!(\sum_{j=1}^n\frac{\partial y_j}{\partial
x_j})\bigwedge_{k=1}^n(dx_k\wedge du_k)\\
\omega_{X_\Sigma}^n & = & n!\bigwedge_{k=1}^n(dx_k\wedge du_k).
\end{eqnarray*}
From this, we see that
\begin{cor}
$\lambda(\mathcal{L}_{[a]})=\lambda(L_h)$ and $h$ is
Hermitian-Einstein if and only the Lagrangian section $L_h$ is
harmonic. In particular, each equivalence class
$[L]\in\mathbbm{L}(Y,W)/\sim$ is represented by a unique harmonic
Lagrangian section.
\end{cor}
On the other hand, the condition for preserving supersymmetry is
given by the following \textit{MMMS equation}, introduced by
Marino-Minasian-Moore-Strominger in \cite{MMMS99} (see also
\cite{LYZ00}):
$$\textrm{Im }e^{\sqrt{-1}\theta}(F_h+\omega_{X_\Sigma})^n=0,$$
for some $\theta\in\mathbb{R}$. Since
$$(F_h+\omega_{X_\Sigma})^n=(\sum_{j=1}^n(dx_j-\sqrt{-1}dy_j(x))\wedge du_j)^n
=\pm(\Omega_Y|_L)\wedge du_1\wedge\ldots\wedge du_n,$$ $h$ satisfies
the MMMS equation with $\theta\in\mathbb{R}$ if and only if $L_h$ is
special Lagrangian with phase $\theta$.

\section{Further remarks}\label{sec5}

We end this paper by several remarks.

\noindent\textbf{1.} For our purposes, we consider $C^2$ hermitian
metrics and Lagrangian sections whose potential are $C^2$ functions.
One can certainly consider metrics and Lagrangians in other
differentiability classes, but then the growth conditions should be
suitably modified.

In particular, when we only require the metrics to be $C^0$,
\textit{singular} Lagrangians can arise as follows. Given a divisor
$\sum_{i=1}^da_iD_i$ in $X_\Sigma$. Then for every $n$-dimensional
cone $\sigma\in\Sigma[n]$, we can find a unique $u_\sigma\in M$ such
that $\langle u_\sigma,v_i\rangle=-a_i$ for all $v_i\in\sigma$. This
defines a piecewise linear function
$\varphi:N_\mathbb{R}\rightarrow\mathbb{R}$ by $\varphi(\xi)=\langle
u_\sigma,\xi\rangle$, for $\xi\in\sigma$. Let $[a]\in
H^2(X_\Sigma,\mathbb{Z})$ be the class represented by
$a=(a_1,\ldots,a_d)$. Then there is a $T_N$-invariant $C^0$
hermitian metric $h$ on the line bundle $\mathcal{L}_{[a]}$ such
that $g_h(\xi)=\varphi(-\xi)$, and $dg_h:N_\mathbb{R}\rightarrow
M_\mathbb{R}$ is the piecewise constant map given by
$dg_h(\xi)=-u_\sigma$ for all $\xi\in\sigma$. Applying the SYZ
mirror transformation, we get a singular Lagrangian
$L_h=\mathcal{F}(h)\subset Y$. This satisfies the following boundary
condition at infinity: for $\xi(t)=t_iv_i+\ldots$, $\langle
dg_h(\xi(t)),v_i\rangle+a_i=0$ for $t_i$ sufficiently negative. In
this case, different line bundles may give rise to the same
Lagrangian subspace. For example, $\mathcal{O}(1)$ and
$\mathcal{O}(-1)$ both transformed to the Lagrangian $L$, which is
the zero section plus the fiber over $\xi=0\in\mathbb{R}$. One can
distinguish the Lagrangian cycles corresponding to $\mathcal{O}(1)$
and $\mathcal{O}(-1)$ by equipping the circle fiber with different
orientations. In this way, the SYZ mirror transformation would still
give a bijective correspondence between isomorphism classes of
holomorphic line bundles over $X_\Sigma$ and equivalence classes of
Lagrangian sections of $(Y,W)$.\\

\noindent\textbf{2.} The SYZ mirror transformation we discuss in
this note only gives a bijective correspondence between holomorphic
line bundles over $X_\Sigma$ and Lagrangian sections of $(Y,W)$. But
it should be extended to an equivalence between the derived category
of coherent sheaves $D^bCoh(X_\Sigma)$ and a suitable variant of the
Fukaya-Kontsevich-Seidel category of $(Y,W)$. In particular, it is
interesting to see how higher rank holomorphic vector bundles over
$X_\Sigma$ can be transformed to Lagrangian multi-sections of
$(Y,W)$ equipped with certain extra data. We plan to address this in
the future.\\

\noindent\textbf{3.} Since we equip $Y$ with the dual of the toric
metric, it is not always possible to represent an equivalence class
$[L]\in\mathbbm{L}(Y,W)/\sim$ by a \textit{minimal} Lagrangian
section. The mirror of $\mathbb{C}P^1$ provides the simplest example
of this. Our way out is to introduce the notion of harmonic
Lagrangians, and as a corollary to our main result, we saw that each
equivalence class $[L]$ is indeed represented by a unique harmonic
representative. However, it would also interesting to look at the
variational theory of Lagrangian sections of $(Y,W)$
\textit{directly}. For example, one may attempt to prove the
existence and uniqueness of harmonic Lagrangian sections by directly
solving the PDE (\ref{harmonic}). On the other hand, the existence
and uniqueness of the solutions of the MMMS equation and the special
Lagrangian equation are largely unexplored. The toric case we
considered here should be the first nontrivial case for one to
investigate these equations.

\appendix
\section{}
In this appendix, we give a proof of Lemma~\ref{lem1}, which is
restated as follows:
\begin{lem}[=Lemma~\ref{lem1}]
$\tilde{L}_{h_0}$ satisfies the growth condition $(\ast_a)$, i.e.
the function $g_{h_0}:N_\mathbb{R}\rightarrow\mathbb{R}$ defined by
$g_{h_0}=-\frac{1}{2}\log h_0(s,s)$ satisfies the following
conditions: Given any $n$-dimensional cone $\sigma\in\Sigma$.
Suppose that, without loss of generality, $\sigma$ is generated by
$v_1,\ldots,v_n$; and let
$\xi(t)=\xi(t_1,\ldots,t_n)=t_1v_1+\ldots+t_nv_n$, for
$t=(t_1,\ldots,t_n)\in\mathbb{R}^n$. Then, we have,
\begin{enumerate}
\item[1.]  the functions
$2e^{-2t_j}(\langle dg(\xi(t)),v_j\rangle+a_j)$ and
$e^{-2t_j}(v_j^T\textrm{Hess}(g)v_j)(\xi(t))$ have the same limit as
$t_j\rightarrow-\infty$, for $j=1,\ldots,n$;
\item[2.] for any $j,k,l\in\{1,\ldots,n\}$, the function
$(v_j^T\textrm{Hess}(g)v_k)(\xi(t))$ has a limit as
$t_l\rightarrow-\infty$; and,
\item[3.] for any distinct $j,k\in\{1,\ldots,n\}$, the function
$e^{-t_j-t_k}(v_j^T\textrm{Hess}(g)v_k)(\xi(t))$ goes to zero when
$t_j\rightarrow-\infty$ or $t_k\rightarrow-\infty$.
\end{enumerate}
\end{lem}
\begin{proof}
By definition, we have $g_{h_0}=-\frac{1}{2}\log
h_0(s,s)=-\frac{1}{2}\sum_{i=1}^da_i\log(l_i\circ\mu_U)$, so that
$$g_{h_0}(\xi)=-\frac{1}{2}\sum_{i=1}^da_i\log\Bigg(\frac{\sum_{u\in\bar{P}\cap
M}c_ul_i(u)e^{2\langle u,\xi\rangle}}{\sum_{u\in\bar{P}\cap
M}c_ue^{2\langle u,\xi\rangle}}\Bigg).$$ The first-order partial
derivatives are given by
$$\frac{\partial g_{h_0}}{\partial\xi_j}=\sum_{i=1}^da_i\Bigg(
\frac{\sum_{u\in\bar{P}\cap M}c_uu^je^{2\langle
u,\xi\rangle}}{\sum_{u\in\bar{P}\cap M}c_ue^{2\langle
u,\xi\rangle}}-\frac{\sum_{u\in\bar{P}\cap M}c_ul_i(u)u^je^{2\langle
u,\xi\rangle}}{\sum_{u\in\bar{P}\cap M}c_ul_i(u)e^{2\langle
u,\xi\rangle}}\Bigg),$$ for $j=1,\ldots,n$. Then, for
$k=1,\ldots,n$,
\begin{eqnarray*}
e^{-2t_k}(\langle dg_{h_0}(\xi(t)),v_k\rangle+a_k) & = &
\sum_{i=1}^da_i\Bigg[\frac{\sum_{l_k(u)\geq1}l_k(u)b_ue^{2(l_k(u)-1)t_k}}{\sum_{u\in\bar{P}\cap
M}b_ue^{2l_k(u)t_k}}\\
&   &
-\frac{\sum_{l_k(u),l_i(u)\geq1}l_k(u)l_i(u)b_ue^{2(l_k(u)-1)t_k}}{\sum_{l_i(u)\geq1}l_i(u)b_ue^{2l_k(u)t_k}}
+\delta_{ik}e^{-2t_k}\Bigg],
\end{eqnarray*}
where
\begin{eqnarray*}
b_u & = & b_u(t_1,\ldots,t_{k-1},t_{k+1},\ldots,t_n)\\
& = &
c_ue^{2(l_1(u)t_1+\ldots+l_{k-1}(u)t_{k-1}+l_{k+1}(u)t_{k+1}+\ldots+l_n(u)t_n)}.
\end{eqnarray*}
Since, for each $i=1,\ldots,d$, there exists $u\in\bar{P}\cap M$
with $l_i(u)=1$ and
$l_1(u)=\ldots=l_{k-1}(u)=l_{k+1}(u)=\ldots=l_n(u)=0$, the limit of
the function $e^{-2t_k}(\langle dg_{h_0}(\xi(t)),v_k\rangle+a_k)$
always exists as $t_m$ goes to $-\infty$ for any
$m=1,\ldots,k-1,k+1,\ldots,n$. Similar arguments apply to other
functions below. Now, as $t_k\rightarrow-\infty$, the terms with the
lowest powers of $e^{t_k}$ dominate. Also note that, for
$i=1,\ldots,d$, there exists $u\in\bar{P}\cap M$ such that
$l_i(u)=1$. So the function $e^{-2t_k}(\langle
dg_{h_0}(\xi(t)),v_k\rangle+a_k)$ has a limit given by
\begin{eqnarray*}
\Big(\sum_{i=1}^da_i\Big)\Bigg(\frac{\sum_{l_k(u)=1}b_u}{\sum_{l_k(u)=0}b_u}\Bigg)
-\sum_{i\neq
k}a_i\Bigg(\frac{\sum_{l_k(u)=1,l_i(u)\geq1}l_i(u)b_u}{\sum_{l_k(u)=0,l_i(u)\geq1}l_i(u)b_u}\Bigg)
-2a_k\Bigg(\frac{\sum_{l_k(u)=2}b_u}{\sum_{l_k(u)=1}b_u}\Bigg).
\end{eqnarray*}

The second-order partial derivatives are given by
\begin{eqnarray*}
\frac{\partial^2g_{h_0}}{\partial\xi_p\partial\xi_q} & = &
2\sum_{i=1}^da_i\Bigg[\frac{\sum_{u\in\bar{P}\cap
M}c_uu^pu^qe^{2\langle u,\xi\rangle}}{\sum_{u\in\bar{P}\cap
M}c_ue^{2\langle u,\xi\rangle}}\\
&   & -\frac{\big(\sum_{u\in\bar{P}\cap M}c_uu^pe^{2\langle
u,\xi\rangle}\big)\big(\sum_{u\in\bar{P}\cap M}c_uu^qe^{2\langle
u,\xi\rangle}\big)}{\big(\sum_{u\in\bar{P}\cap M}c_ue^{2\langle
u,\xi\rangle}\big)^2}\\
&   & -\frac{\sum_{u\in\bar{P}\cap M}c_ul_i(u)u^pu^qe^{2\langle
u,\xi\rangle}}{\sum_{u\in\bar{P}\cap M}c_ul_i(u)e^{2\langle
u,\xi\rangle}}\\
&   & +\frac{\big(\sum_{u\in\bar{P}\cap M}c_ul_i(u)u^pe^{2\langle
u,\xi\rangle}\big)\big(\sum_{u\in\bar{P}\cap
M}c_ul_i(u)u^qe^{2\langle
u,\xi\rangle}\big)}{\big(\sum_{u\in\bar{P}\cap
M}c_ul_i(u)e^{2\langle u,\xi\rangle}\big)^2}\Bigg],
\end{eqnarray*}
for $p,q=1,\ldots,n$. From this we compute, for
$j,k\in\{1,\ldots,n\}$,
\begin{eqnarray*}
&   & v_j^T\textrm{Hess}(g_{h_0})v_k\\
& = & \sum_{p,q=1}^nv_j^pv_k^q\frac{\partial^2g_{h_0}}{\partial\xi_p\partial\xi_q}\\
& = &
2\sum_{i=1}^da_i\Bigg[\frac{\sum_{l_j(u),l_k(u)\geq1}c_ul_j(u)l_k(u)e^{2\langle
u,\xi\rangle}}{\sum_{u\in\bar{P}\cap M}c_ue^{2\langle
u,\xi\rangle}}\\
&   & -\frac{\big(\sum_{l_j(u)\geq1}c_ul_j(u)e^{2\langle
u,\xi\rangle}\big)\big(\sum_{l_k(u)\geq1}c_ul_k(u)e^{2\langle
u,\xi\rangle}\big)}{\big(\sum_{u\in\bar{P}\cap M}c_ue^{2\langle
u,\xi\rangle}\big)^2}\\
&   &
-\frac{\sum_{l_i(u),l_j(u),l_k(u)\geq1}c_ul_i(u)l_j(u)l_k(u)e^{2\langle
u,\xi\rangle}}{\sum_{l_i(u)\geq1}c_ul_i(u)e^{2\langle
u,\xi\rangle}}\\
&   &
+\frac{\big(\sum_{l_i(u),l_j(u)\geq1}c_ul_i(u)l_j(u)e^{2\langle
u,\xi\rangle}\big)\big(\sum_{l_i(u),l_k(u)\geq1}c_ul_i(u)l_k(u)e^{2\langle
u,\xi\rangle}\big)}{\big(\sum_{l_i(u)\geq1}c_ul_i(u)e^{2\langle
u,\xi\rangle}\big)^2}\Bigg].
\end{eqnarray*}
It is easy to see that, for any $l=1,\ldots,n$, as
$t_l\rightarrow-\infty$, the limit of the function
$(v_j^T\textrm{Hess}(g_{h_0})v_k)(\xi(t))$ always exists. When
$j\neq k$, let
\begin{eqnarray*}
b_u & = &
b_u(t_1,\ldots,\widehat{t}_j,\ldots,\widehat{t}_k,\ldots,t_n)\\
& = &
c_ue^{2(l_1(u)t_1+\ldots+\widehat{l_j(u)t_j}+\ldots+\widehat{l_k(u)t_k}+\ldots+l_n(u)t_n)}.
\end{eqnarray*}
Then the function $(v_j^T\textrm{Hess}(g_{h_0})v_k)(\xi(t))$ is
equal to the following expression
\begin{eqnarray*}
&   &
2\sum_{i=1}^da_i\Bigg[\frac{\sum_{l_j(u),l_k(u)\geq1}l_j(u)l_k(u)b_ue^{2(l_j(u)t_j+l_k(u)t_k)}}
{\sum_{u\in\bar{P}\cap M}b_ue^{2(l_j(u)t_j+l_k(u)t_k)}}\\
&   &
-\frac{\big(\sum_{l_j(u)\geq1}l_j(u)b_ue^{2(l_j(u)t_j+l_k(u)t_k)}\big)
\big(\sum_{l_k(u)\geq1}l_k(u)b_ue^{2(l_j(u)t_j+l_k(u)t_k)}\big)}
{\big(\sum_{u\in\bar{P}\cap M}b_ue^{2(l_j(u)t_j+l_k(u)t_k)}\big)^2}\\
&   &
-\frac{\sum_{l_i(u),l_j(u),l_k(u)\geq1}l_i(u)l_j(u)l_k(u)b_ue^{2(l_j(u)t_j+l_k(u)t_k)}}
{\sum_{l_i(u)\geq1}l_i(u)b_ue^{2(l_j(u)t_j+l_k(u)t_k)}}\\
&   &
+\frac{\big(\sum_{l_i(u),l_j(u)\geq1}l_i(u)l_j(u)b_ue^{2(l_j(u)t_j+l_k(u)t_k)}\big)
\big(\sum_{l_i(u),l_k(u)\geq1}l_i(u)l_k(u)b_ue^{2(l_j(u)t_j+l_k(u)t_k)}\big)}
{\big(\sum_{l_i(u)\geq1}l_i(u)b_ue^{2(l_j(u)t_j+l_k(u)t_k)}\big)^2}\Bigg].
\end{eqnarray*}
Notice that in each term, the numerator is $O(e^{2t_j+2t_k})$, while
the denominator is $O(1)$. Thus, the function
$e^{-t_j-t_k}(v_j^T\textrm{Hess}(g_{h_0})v_k)(\xi(t))$ goes to zero
as $t_j\rightarrow-\infty$ or $t_k\rightarrow-\infty$.

For $j=k$, let
\begin{eqnarray*}
b_u & = & b_u(t_1,\ldots,t_{k-1},t_{k+1},\ldots,t_n)\\
& = &
c_ue^{2(l_1(u)t_1+\ldots+l_{k-1}(u)t_{k-1}+l_{k+1}(u)t_{k+1}+\ldots+l_n(u)t_n)}.
\end{eqnarray*}
Then
\begin{eqnarray*}
e^{-2t_k}(v_k^T\textrm{Hess}(g_{h_0})v_k)(\xi(t)) & = &
2e^{-2t_k}\sum_{i=1}^da_i\Bigg[\frac{\sum_{l_k(u)\geq1}l_k(u)^2b_ue^{2l_k(u)t_k}}{\sum_{u\in\bar{P}\cap
M}b_ue^{2l_k(u)t_k}}\\
&   &
-\Bigg(\frac{\sum_{l_k(u)\geq1}l_k(u)b_ue^{2l_k(u)t_k}}{\sum_{u\in\bar{P}\cap
M}b_ue^{2l_k(u)t_k}}\Bigg)^2\\
&   &
-\frac{\sum_{l_i(u),l_k(u)\geq1}l_i(u)l_k(u)^2b_ue^{2l_k(u)t_k}}{\sum_{l_i(u)\geq1}l_i(u)b_ue^{2l_k(u)t_k}}\\
&   &
+\Bigg(\frac{\sum_{l_i(u),l_k(u)\geq1}l_i(u)l_k(u)b_ue^{2l_k(u)t_k}}
{\sum_{l_i(u)\geq1}l_i(u)b_ue^{2l_k(u)t_k}}\Bigg)^2\Bigg].
\end{eqnarray*}
As $t_k\rightarrow-\infty$, the function
$e^{-2t_k}(v_k^T\textrm{Hess}(g_{h_0})v_k)(\xi(t))$ has a limit
given by
\begin{eqnarray*}
2\Big(\sum_{i=1}^da_i\Big)\Bigg(\frac{\sum_{l_k(u)=1}b_u}{\sum_{l_k(u)=0}b_u}\Bigg)
-2\sum_{i\neq
k}a_i\Bigg(\frac{\sum_{l_k(u)=1,l_i(u)\geq1}l_i(u)b_u}{\sum_{l_k(u)=0,l_i(u)\geq1}l_i(u)b_u}\Bigg)
-4a_k\Bigg(\frac{\sum_{l_k(u)=2}b_u}{\sum_{l_k(u)=1}b_u}\Bigg).
\end{eqnarray*}
This completes the proof of Lemma~\ref{lem1}.
\end{proof}


\begin{thebibliography}{99}

\bibitem{Abouzaid05}
M. Abouzaid, \emph{Homogeneous coordinate rings and mirror symmetry
for toric varieties}. Geom. Topol., {\bf 10} (2006), 1097--1157
(math.SG/0511644).

\bibitem{Abouzaid06}
\underline{\qquad\quad}, \emph{Morse homology, tropical geometry,
and homological mirror symmetry for toric varieties}. Preprint, 2006
(math/0610004).

\bibitem{Auroux07}
D. Auroux, \emph{Mirror symmetry and T-duality in the complement of
an anticanonical divisor}. J. Gokova Geom. Topol. GGT, {\bf 1}
(2007), 51--91 (arXiv:0706.3207).

\bibitem{AKO04}
D. Auroux, L. Katzarkov and D. Orlov, \emph{Mirror symmetry for
weighted projective planes and their noncommutative deformations}.
Ann. of Math. (2), {\bf 167} (2008), no. 3, 867--943
(math.AG/0404281).

\bibitem{AKO05}
\underline{\qquad\quad}, \emph{Mirror symmetry for del Pezzo
surfaces: vanishing cycles and coherent sheaves}. Invent. Math.,
{\bf 166} (2006), no. 3, 537--582 (math.AG/0506166).

\bibitem{Chan-Leung08a}
K.-W. Chan and N.-C. Leung, \emph{Mirror symmetry for toric Fano
manifolds via SYZ transformations}. Preprint, 2008
(arXiv:0801.2830).

\bibitem{Chan-Leung08b}
\underline{\qquad\quad}, \emph{On SYZ mirror transformations}. To
appear in Advanced Studies in Pure Mathematics, "New developments in
Algebraic Geometry, Integrable Systems and Mirror Symmetry".
(arXiv:0808.1551).

\bibitem{Fang08}
B. Fang, \emph{Homological mirror symmetry is T-duality for
$\mathbb{P}^n$}. Commun. Number Theory Phys. 2 (2008), no. 4,
719--742 (arXiv:0804.0646).

\bibitem{FLTZ08}
B. Fang, C.-C. M. Liu, D. Treumann and E. Zaslow, \emph{A
categorification of Morelli's theorem and homological mirror
symmetry for toric varieties}. Preprint, 2008 (arXiv:0811.1228).

\bibitem{Fulton93}
W. Fulton, \emph{Introduction to toric varieties}. Annals of
Mathematics Studies, 131. The William H. Roever Lectures in
Geometry. Princeton University Press, Princeton, NJ, 1993.

\bibitem{Guillemin94a}
V. Guillemin, \emph{Moment maps and combinatorial invariants of
Hamiltonian $T^n$-spaces}. Progress in Mathematics, 122.
Birkha\"{u}ser Boston, Inc., Boston, MA, 1994.

\bibitem{Guillemin94b}
\underline{\qquad\quad}, \emph{Kaehler structures on toric
varieties}. J. Differential Geom. {\bf 40} (1994), no. 2, 285--309.


\bibitem{HIV00}
K. Hori, A. Iqbal and C. Vafa, \emph{D-branes and mirror symmetry}.
Preprint, 2000 (hep-th/0005247).

\bibitem{Kontsevich98}
M. Kontsevich, \emph{Lectures at ENS}. Paris, Spring 1998, notes
taken by J. Bellaiche, J.-F. Dat, I. Marin, G. Racinet and H.
Randriambololona.

\bibitem{LYZ00}
N.-C. Leung, S.-T. Yau and E. Zaslow, \emph{From special Lagrangian
to Hermitian-Yang-Mills via Fourier-Mukai transform}. Adv. Theor.
Math. Phys., {\bf 4} (2000), no. 6, 1319--1341 (math.DG/0005118).

\bibitem{MMMS99}
M. Mari\~{n}o, R. Minasian, G. Moore and A. Strominger,
\emph{Nonlinear instantons from supersymmetric $p$-branes}. JHEP 01
(2000), 005 (hep-th/9911206).

\bibitem{Seidel00}
P. Seidel, \emph{More about vanishing cycles and mutation}.
Symplectic geometry and mirror symmetry (Seoul, 2000), 429--465,
World Sci. Publ., River Edge, NJ, 2001 (math.SG/0010032).

\bibitem{SYZ96}
A. Strominger, S.-T. Yau and E. Zaslow, \emph{Mirror symmetry is
T-duality}. Nuclear Phys. B, {\bf 479} (1996), no. 1-2, 243--259
(hep-th/9606040).

\bibitem{Ueda04}
K. Ueda, \emph{Homological mirror symmetry for toric del Pezzo
surfaces}. Comm. Math. Phys., {\bf 264} (2006), no. 1, 71--85
(math.AG/0411654).

\end{thebibliography}
\end{document}